\documentclass[12pt,english]{article}
\usepackage{amsthm,amsfonts,amsmath,amssymb}
\usepackage[utf8]{inputenc}
\usepackage[english, russian]{babel}
\usepackage{srcltx}
\usepackage{latexsym,,longtable,mathrsfs,array}
\usepackage[unicode,final,hyperindex,hypertex]{hyperref}
\usepackage{bbm}
\usepackage{lineno}

\usepackage{color}

\usepackage{mathabx}
\usepackage{enumerate}
\usepackage[centerlast, small]{caption}

\newtheorem{thm}{Theorem}
\newtheorem{lem}{Lemma}

\newcommand{\F}{\mathbbmss{F}}    
\newcommand{\splitext}{\,\colon\!}
\newcommand{\arbitraryext}{\,\ldotp}
\newcommand{\GL}{{\mathrm{GL}}}

\newcommand{\PSL}{\mathop{\rm PSL}\nolimits}           
\newcommand{\SL}{\mathop{\rm SL}\nolimits}           
\newcommand{\SU}{\mathop{\rm SU}\nolimits}             
\newcommand{\PSU}{\mathop{\rm PSU}\nolimits}  
\newcommand{\PGU}{\mathop{\rm PGU}\nolimits}
\newcommand{\PGL}{{\mathrm{PGL}}} 
\renewcommand{\P}{{\rm P}}

\newcommand{\Aut}{\operatorname{Aut}}

\newcommand{\ld}{\mathop{\ldbrack}}
\newcommand{\rd}{\mathop{\rdbrack}}

\newcommand{\D}{\mathcal{D}}
\renewcommand{\C}{\mathcal{C}}
\newcommand{\E}{\mathcal{E}}
\renewcommand{\U}{\mathcal{U}}

\newcommand{\lcm}{\mathrm{lcm}}

\begin{document}
 \selectlanguage{english}

\title{ On the heritability of the Hall property $\D_\pi$ by overgroups of $\pi$-Hall subgroups \footnote{The work is supported by  Russian Science Foundation (project 14-21-00065).}}

\author{N.Ch. Manzaeva}

\date{}
\maketitle

\begin{abstract}
 Let $\pi$ be a set of primes. We say that a finite group $G$ is a $\D_\pi$-group if all its maximal $\pi$-subgroups are conjugate. In this paper, we give an affirmative answer to the problem 17.44(b) from ``The Kourovka notebook'': we prove that in a $\D_\pi$-group an overgroup of a $\pi$-Hall subgroup is always a $\D_\pi$-group.
\end{abstract}

\section{Introduction}

Let $G$~be a finite group, $\pi$~ a set of primes.  We denote by $\pi'$ the set of all primes not in $\pi$, by $\pi(n)$  the set of all prime divisors of a positive integer~$n$, given a group $G$ we denote $\pi(|G|)$ by $\pi(G)$. A group $G$ with $\pi(G)\subseteq \pi$ is called a {\it $\pi$-group}. A subgroup $H$ of $G$ is called a {\it $\pi$-Hall subgroup}, if $\pi(H)\subseteq \pi$  and $\pi(|G:H|)\subseteq \pi'$.

Following \cite{Hall1956}, we say that $G$
{\it satisfies $\E_\pi$} (or briefly $G\in \E_\pi$), if $G$ has a $\pi$-Hall subgroup. If $G$ satisfies $\E_\pi$ and every two $\pi$-Hall subgroups of $G$ are conjugate, then we say that $G$ {\it satisfies $\C_\pi$} ($G\in \C_\pi$). Finally, $G$ {\it satisfies $\D_\pi$} ($G\in \D_\pi$), if $G$ satisfies $\C_\pi$ and every $\pi$-subgroup of $G$ is included in a $\pi$-Hall subgroup of~$G$. Thus $G\in \D_\pi$ if a complete analogue of the Sylow theorem for $\pi$-subgroups of $G$ holds. A group satisfying $\E_\pi$ ($\C_\pi$, $\D_\pi$) is also called an {\it$\E_\pi$-group}
(respectively $\C_\pi$-group, $\D_\pi$-group).

 In the paper,  we give an affirmative answer to the following problem from ``The Kourovka notebook'' \cite{Kour}:

\bigskip

\noindent{\bfseries Problem 1.}\label{17.44b}{\cite[Problem 17.44(b)]{Kour}}
 In a $\D_\pi$-group, is an overgroup of a $\pi$-Hall subgroup always a $\D_\pi$-group?

\bigskip

An affirmative answer to the analogous problem for $\C_\pi$-property (see \cite[Problem 17.44(a)]{Kour}) is obtained by E.P.~Vdovin and D.O. Revin in \cite{VdovinRevin2012, VdovinRevin2013}.

According to \cite{VdovinRevin2011}, we say that $G$ {\it satisfies $\U_\pi$}, if $G\in \D_\pi$ and every overgroup of a $\pi$-Hall subgroup of $G$ satisfies $\D_\pi$. Thus Problem~\ref{17.44b} can be reformulated in the following way:

\bigskip

\noindent{\bfseries Problem 1.}  Is it true that $\D_\pi=\U_\pi$?

\bigskip

The following main theorem gives an affirmative answer to Problem~1.

\begin{thm}{\em(Main theorem)}\label{main}
Let $\pi$~be a set of primes.  Then  $\D_\pi=\U_\pi$. In other words, if $G\in \D_\pi$ and $H$ is a $\pi$-Hall subgroup of $G$, then every subgroup~$M$ of $G$ with $H\le M$ satisfies~$\D_\pi$.
\end{thm}

 In \cite[Theorem~7.7]{RevinVdovin2006}, E.P.~Vdovin and D.O.~Revin  proved that  $G$ satisfies~$\D_\pi$ if and only if each composition factor of $G$ satisfies $\D_\pi$. Using this result, an analogous criterion for $\U_\pi$ is obtained in \cite{VdovinManRevin2012}.

 \begin{thm}\label{reduction}{\em \cite[Theorem~2]{VdovinManRevin2012}}
A finite group $G$ satisfies $\U_\pi$ if and only if each composition factor of $G$ satisfies~$\U_\pi$.
 \end{thm}

Thus Problem 1 is reduced to a similar problem for simple $\D_\pi$-groups. All simple $\D_\pi$-groups are known: pure arithmetic necessary and sufficient conditions for a simple group $G$ to satisfy $\D_\pi$ can be found in \cite{Revin2008e}. It was proved in   \cite{VdovinManRevin2012} that if $G\in \D_\pi$ is an alternating group, a sporadic group or a group of Lie type in characteristic $p\in \pi$, then $G$ satisfies~$\U_\pi$. An affirmative answer to Problem~\ref{17.44b} in case $2\in \pi$ is obtained in \cite{Manzaeva2014}. Therefore, in this paper, we consider only   $\D_\pi$-groups of Lie type in characteristic~$p$ with $2,p\notin \pi$.

\section{Notation and preliminary results}

All groups in the paper are assumed to be finite. Our notation is standard and agrees with that of \cite{CFSG} and \cite{KleiLie}. By $A\splitext B$ and $A\arbitraryext B$ we denote a split extension and an arbitrary extension of a group $A$ by a group $B$, respectively. Symbols $A\times B$ and $A\circ B$ denote direct and central products, respectively. If $G$ is a group and $S$ is a permutation group, then $G\wr S$ is the permutation wreath product of $G$ and~$S$. For $M\subseteq G$ we set $M^G=\{M^g\mid g\in G\}$. The subgroup generated by a subset $M$ is denoted by $\langle M \rangle$. We use notations $H\leq G$ and $H \unlhd G$ instead of ``$H$ is a subgroup of~$G$'' and ``$H$ is a normal subgroup of~$G$'', respectively.   The normalizer and the centralizer of $H$ in $G$ are denoted by $N_G(H)$ and $C_G(H)$, respectively, while $Z(G)$ is the center of~$G$.  The generalized Fitting subgroup of $G$ is denoted by $F^*(G)$. Denote by $\ldbrack n \rdbrack$ an arbitrary solvable group of order $n$.

Throughout, $\F_q$ is a finite field of order $q$ and characteristic $p$. By $\eta$ we always denote an element from the set $\{+, -\}$ and we use $\eta$ instead of $\eta1$ as well. In order to unify notation and arguments we often denote $A_n(q)$ by $A_n^+(q)$, ${}^2A_n(q)$ by $A_n^-(q)$, $E_6(q)$ by $E_6^+(q)$, ${}^2E_6(q)$ by $E_6^-(q)$, $\PSL_n(q)$ by $\PSL_n^+(q)$ and $\PSU_n(q)$ by $\PSL_n^-(q)$.  If $G$ is a group of Lie type, then by $W(G)$ we denote the Weyl group of $G$.

 The integral part of a real number $x$ is denoted by   $\left[ x \right]$. For integers $n$ and  $m$, we denote by $\gcd(n,m)$ and $\lcm(n,m)$ the greatest common divisor and the least common multiple, respectively. If  $n$~is a positive integer, then $n_\pi$ is the largest divisor $d$ of $n$ with $\pi(d)\subseteq\pi$.  If $r$~is an odd prime and  $q$~is an integer not divisible by $r$, then $e(q,r)$ is the smallest positive integer  $e$ with $q^e \equiv 1 \pmod r$.

We say that $G$ {\it satisfies} $(*)$ if every $\pi$-subgroup of $G$ has a normal abelian $\tau$-Hall subgroup, where $\tau=(\pi\cap\pi(G))\setminus \{r\}$ and $r=min(\pi\cap\pi(G))$.

The next result may be found in \cite{Weir1955}.

\begin{lem}\label{r4ast'}{\em ( \cite{Weir1955}, \cite[Lemmas~2.4 and 2.5]{Gross1995})}
 Let $r$ be an odd prime, $k$~an integer not divisible by $r$,  and  $m$~a positive integer. Denote  $e(k,r)$ by $e$ and

 $$e^*=\left\{
\begin{array}{ll}
2e & \text{ if }  e\equiv 1\pmod 2,\\
e & \mbox{ if   }  e\equiv 0\pmod 4,\\
e/2 & \mbox{ if  }  e\equiv 2\pmod 4.
\end{array}
\right.
$$

 Then the following identities hold:

$$(k^m-1)_r=\left\{
\begin{array}{ll}
(k^e-1)_r (m/e)_{r} & \text{ if } e \text{ divides } m,\\
1 & \text{ if } e \text{ does not divide } m;\end{array}\right.$$

$$(k^m-(-1)^m)_r=\left\{
\begin{array}{ll}
(k^{e^*}-(-1)^{e^*})_r (m/{e^*})_{r} & \text{ if } e^* \text{ divides } m,\\
1 & \text{ if }  e^* \text{ does not divide } m.\end{array}\right.$$


\end{lem}

In Lemma~\ref{L_Epi} we collect some known facts about $\pi$-Hall subgroups in finite groups.

\begin{lem}\label{L_Epi}
Let $G$ be a finite group, $A$~a normal subgroup of $G$.
\begin{enumerate}[{\em(a)}]
 \item If $H$ is a $\pi$-Hall subgroup of $G$, then $H\cap A$ is a $\pi$-Hall subgroup of $A$ and $HA/A$~is a $\pi$-Hall subgroup of $G/A$. In particular, a normal subgroup and a homomorphic image of an  $\E_\pi$-group satisfy $\E_\pi$. {\em  (see \cite[Lemma~1]{Hall1956})} \label{Epi}


\item If $M/A$ is a $\pi$-subgroup of $G/A$, then there exists a $\pi$-subgroup $H$ of $G$ with $M=HA$. {\em (see  \cite[Lemma~2.1]{Gross1986})} \label{existpi}

\item If $M$ is~a $\pi$-Hall subgroup of $A$ and $G/A$ is a $\pi$-group, then a $\pi$-Hall subgroup $H$ of $G$ with $H\cap A=M$ exists if and only if the class  $M^A$ is $G$-invariant, i.e., ${M^A=M^G}$. {\em (see  \cite[Proposition~4.8]{VdovinRevin2011})} \label{prop4.8}

\item If $2\notin \pi$ then $\E_\pi=\C_\pi$. In particular, a group $G$ satisfies $\E_\pi$ if and only if each composition factor of $G$ satisfies $\E_\pi$. {\em  (see \cite[Theorem A]{Gross1987}, \cite[Theorem 2.3]{Gross1995}, \cite[Theorem~5.4]{VdovinRevin2011})} \label{epi=cpi}

\item If $G$ possesses a nilpotent $\pi$-Hall subgroup then $G$ satisfies $\D_\pi$. {\em (see  \cite{Wielandt1954}, \cite[Theorem~6.2]{VdovinRevin2011}) } \label{dpi_nil}

\item A group~$G$ satisfies $\D_\pi$ if and only if  $A$ and $G/A$ satisfy $\D_\pi$. Equivalently,  $G\in \D_\pi$ if and only if each composition factor of $G$ satisfies~$\D_\pi$. {\em (see   \cite[Theorem~7.7]{RevinVdovin2006}, \cite[Collorary~6.7]{VdovinRevin2011})} \label{Dpi_comp_crit}

\end{enumerate}
\end{lem}

\begin{lem} \label{Dpi}
{\em (\cite[Theorem~3]{Revin2008e}, \cite[Theorem~6.9]{VdovinRevin2011})}
Let $S$ be a simple group of Lie type with the base field  $\F_q$ of characteristic $p$. Suppose $2\not \in \pi$ and $|\pi\cap\pi(S)|\ge 2$. Then  $S$ satisfies~$\D_\pi$ if and only if the pair~$(S, \pi)$ satisfies one of the Condition  I-IV below.
\end{lem}

\noindent{\bf Condition I.} Let  $p\in\pi$ and $\tau=(\pi\cap\pi(S))\setminus\{p\}$. We say that~$(S,\pi)$ {\it satisfies Condition}~I if
$\tau\subseteq\pi(q-1)$ and every number from $\pi$ does not divide~$|W(S)|$.

\smallskip

\noindent{\bf Condition II.} Suppose that $S$ is not isomorphic to ${^2B}_2(q),{^2F}_4(q),{^2G}_2(q)$ and $p\not\in\pi$. Denote by  $r$ the number $min(\pi\cap\pi(S))$. Set  $\tau=(\pi\cap\pi(S))\setminus\{r\}$ and
$a=e(q,r)$. We say that $(S,\pi)$ {\it satisfies Condition}~II if there exists  $t\in\tau$ with
$b=e(q,t)\ne a$ and one of the following holds.
\begin{enumerate}[(a)]

\item  $S\simeq A_{n-1}(q)$, $a=r-1$, $b=r$, $(q^{r-1}-1)_r=r$,
$\left[\displaystyle\frac{n}{r-1}\right]=\left[\displaystyle\frac{n}{r}\right]$
and for every $s\in\tau$ hold $e(q,s)=b$ and $n<bs$.

\item  $S\simeq A_{n-1}(q)$, $a=r-1$, $b=r$, $(q^{r-1}-1)_r=r$,
$\left[\displaystyle\frac{n}{r-1}\right]=\left[\displaystyle\frac{n}{r}\right]+1$,
$n\equiv -1 \pmod r$ and for every $s\in\tau$ hold $e(q,s)=b$ and $n<bs$.

\item  $S\simeq {^2A}_{n-1}(q)$, $r\equiv 1 \pmod 4$, $a=r-1$, $b=2r$, $(q^{r-1}-1)_r=r$,
$\left[\displaystyle\frac{n}{r-1}\right]=\left[\displaystyle\frac{n}{r}\right]$
and $e(q,s)=b$ for every $s\in\tau$.

\item  $S\simeq {^2A}_{n-1}(q)$, $r\equiv 3 \pmod 4$, $a=\displaystyle\frac{r-1}{2}$,
$b=2r$, $(q^{r-1}-1)_r=r$,
$\left[\displaystyle\frac{n}{r-1}\right]=\left[\displaystyle\frac{n}{r}\right]$
and $e(q,s)=b$ for every $s\in\tau$.

\item  $S\simeq {^2A}_{n-1}(q)$, $r\equiv 1 \pmod 4$, $a=r-1$, $b=2r$, $(q^{r-1}-1)_r=r$,
$\left[\displaystyle\frac{n}{r-1}\right]=\left[\displaystyle\frac{n}{r}\right]+1$,
$n\equiv -1 \pmod r$ and
$e(q,s)=b$ for every $s\in\tau$.

\item $S\simeq {^2A}_{n-1}(q)$, $r\equiv 3 \pmod 4$, $a=\displaystyle\frac{r-1}{2}$,
$b=2r$, $(q^{r-1}-1)_r=r$,
$\left[\displaystyle\frac{n}{r-1}\right]=\left[\displaystyle\frac{n}{r}\right]+1$,
$n\equiv -1 \pmod r$ and
$e(q,s)=b$ for every $s\in\tau$.

\item $S\simeq {^2D}_n(q)$, $a\equiv 1 \pmod 2$, $n=b=2a$ and for every $s\in\tau$ either $e(q,s)=a$ or
$e(q,s)=b$. \label{CondII_2D_n_1}

\item $S\simeq {^2D}_n(q)$, $b\equiv 1 \pmod 2$, $n=a=2b$ and for every $s\in\tau$ either $e(q,s)=a$ or
$e(q,s)=b$. \label{CondII_2D_n_2}

\end{enumerate}

In cases (\ref{CondII_2D_n_1})-(\ref{CondII_2D_n_2}), a $\pi$-Hall subgroup of $S\simeq {}^2D_n(q) $ is cyclic.

\smallskip

\noindent{\bf Condition III.} Suppose that $S$ is not isomorphic to ${^2B}_2(q),{^2F}_4(q),{^2G}_2(q)$ and $p\not\in\pi$. Denote by  $r$ the number $min(\pi\cap\pi(S))$. Set  $\tau=(\pi\cap\pi(S))\setminus\{r\}$ and
$c=e(q,r)$. We say that $(S,\pi)$ {\it satisfies Condition}~III if ${e(q,t)=c}$ for every $t\in\tau$
 and one of the following holds.
\begin{enumerate}[(a)]
\item  $S\simeq A_{n-1}(q)$ and $n<cs$ for every $s\in\tau$.

\item  $S\simeq{^2A}_{n-1}(q)$, $c\equiv 0 \pmod 4$ and $n<cs$ for every $s\in\tau$.

\item  $S\simeq{^2A}_{n-1}(q)$, $c\equiv 2 \pmod 4$ and $2n<cs$ for every $s\in\tau$.

\item  $S\simeq{^2A}_{n-1}(q)$, $c\equiv 1 \pmod 2$ and $n<2cs$ for every $s\in\tau$.

\item  $S$ is isomorphic to one of the groups $B_n(q)$, $C_n(q)$ or ${^2D}_n(q)$,  $c$ is even and $2n<cs$ for every $s\in\tau$.

\item   $S$ is isomorphic to one of the groups $B_n(q)$, $C_n(q)$ or $D_n(q)$, $c$ is odd and $n<cs$ for every $s\in\tau$.

\item  $S\simeq D_n(q)$,  $c$ is even and $2n\le cs$ for every $s\in\tau$.

\item  $S\simeq {^2D}_n(q)$,  $c$  is odd and $n\le cs$ for every $s\in\tau$.

\item  $S\simeq {^3D}_4(q)$.

\item  $S\simeq E_6(q)$  and  if $r=3$ and $c=1$ then $5,13\not\in\tau$.

\item  $S\simeq {^2E}_6(q)$ and if $r=3$ and  $c=2$ then $5,13\not\in\tau$.

\item  $S\simeq E_7(q)$ and if $r=3$ and $c\in\{1,2\}$
then $5,7,13\not\in\tau$, and if $r=5$ and $c\in\{1,2\}$ then
$7\not\in\tau$.

\item  $S\simeq E_8(q)$ and if $r=3$ and $c\in\{1,2\}$ then $5,7,13\not\in\tau$, and if $r=5$ and
$c\in\{1,2\}$ then $7,31\not\in\tau$.

\item  $S\simeq G_2(q)$.

\item  $S\simeq F_4(q)$ and if $r=3$ and $c=1$ then $13\not\in\tau$.

\end{enumerate}

\smallskip

\noindent{\bf Condition IV.} We say that $(S,\pi)$ {\it satisfies Condition}~IV if one of the following holds.

\begin{enumerate}[(a)]

\item $S\simeq {^2B}_2(2^{2m+1})$, $\pi\cap\pi(G)$ is contained in one of the sets
$\pi(2^{2m+1}-1)$, $\pi(2^{2m+1}\pm 2^{m+1}+1)$.

\item $S\simeq {^2G}_2(3^{2m+1})$, $\pi\cap\pi(G)$ is contained in one of the sets
$\pi(3^{2m+1}-1)\setminus\{2\}$, $\pi(3^{2m+1}\pm
3^{m+1}+1)\setminus\{2\}$.

\item $S\simeq {^2F}_4(2^{2m+1})$, $\pi\cap\pi(G)$ is contained in one of the sets
$\pi(2^{2(2m+1)}\pm 1)$, $\pi(2^{2m+1}\pm 2^{m+1}+1)$,
$\pi(2^{2(2m+1)}\pm 2^{3m+2}\mp2^{m+1}-1)$, $\pi(2^{2(2m+1)}\pm
2^{3m+2}+2^{2m+1}\pm2^{m+1}-1)$.

\end{enumerate}

In the next lemmas, we recall some preliminary results about $\U_\pi$-property.

\begin{lem}\label{max}{\em \cite[Lemma~3]{Manzaeva2014}}
The following statements are equivalent.
\begin{enumerate}[{\em(a)}]
\item $\D_\pi=\U_\pi$.
\item In every simple $\D_\pi$-group $G$ all maximal subgroups containing a  $\pi$-Hall subgroup of $G$ satisfy $\D_\pi$.
\end{enumerate}
\end{lem}

\begin{lem}\label{pinpi}{\em \cite[Theorem 4]{VdovinManRevin2012}}
If $G\in\D_\pi$~is either an alternating group or a sporadic simple group or a simple group of Lie type in characteristic $p \in \pi$, then $G$ satisfies $\U_\pi$.
\end{lem}

\begin{lem}\label{2inpi}{\em \cite[Theorem 1]{Manzaeva2014}}
If $2\in \pi$ then  $\D_\pi= \U_\pi$.
\end{lem}

 In view of Lemma~\ref{2inpi}, we consider the case where $\pi$ is a set of odd primes. The next results are concerned with the properties $\E_\pi$ and $\D_\pi$, where $2\notin \pi$.

\begin{lem} \label{t.2.1}  {\em \cite[Theorem 1]{VdovinRevin2002}}
  Let $G$~ be a group of Lie type in characteristic~$p$. Suppose that $2,p\notin\pi$ and $H$~is a~$\pi$-Hall subgroup of~$G$. Set $r= min(\pi\cap\pi(G))$ and  $\tau=\pi\setminus\{r\}$. Then $H$ has a normal abelian  $\tau$-Hall subgroup.
\end{lem}

Recall that $G$  satisfies~$(*)$ if every $\pi$-subgroup of $G$ has a normal abelian $\tau$-Hall subgroup, where $\tau=(\pi\cap\pi(G))\setminus \{r\}$ and $r=min(\pi\cap\pi(G))$. Suppose that $G\in \E_\pi$ is a simple group of Lie type in characteristic~$p$. Lemma~\ref{t.2.1} implies that if $2,p \notin \pi$ and $G \in \D_\pi$ then $G$ satisfies $(*)$. If $e(q,s)=e(q,r)$ for every $s\in \tau$, then the converse is also true.

\begin{lem}  \label{dpi_not2}{\em \cite[Theorem 5]{VdovinRevin2002}}
Let $G$~be a simple group of Lie type with the base field $\F_q$ of characteristic $p$, not isomorphic to ${^2B}_2(q),{^2F}_4(q),{^2G}_2(q)$. Suppose that $2,p \not \in \pi$ and  $G\in \E_\pi$. Set $r=min(\pi\cap \pi(G))$ and $\tau=(\pi\cap \pi(G))\setminus\{r\}$. Assume that $e(q,s)=e(q,r)$ for every $s\in\tau$. Then $G\in \D_\pi$ if and only if $G$ satisfies $(*)$.
\end{lem}

\begin{lem}\label{epi-dpi_lie} {\em (\cite[Theorem~1.1]{Gross1986}, \cite[Theorem~6.14]{Gross1986}, \cite[Lemmas~5-7]{Revin2009around})}
Assume that $2\notin \pi$, $G$ is a simple group and $G \in \E_\pi\setminus \D_\pi$. Set $r=min(\pi\cap \pi(G))$ and $\tau=(\pi\cap \pi(G))\setminus\{r\}$. Then one of the following holds.

\begin{enumerate}[{\em(1)}]
\item $G=O'N$ and $\pi\cap\pi(G)=\{3,5\}$.
\item $G$~is a group of Lie type with the base field $\F_q$ of characteristic $p>0$ and either {\em(A)} or {\em(B)} is true:
    \begin{enumerate}[{\em(A)}]
    \item $p\in \pi$, $p$ divides $|W(G)|$, every $t\in(\pi\cap \pi(G))\setminus \{p\}$ divides $q-1$ and does not divide $|W(G)|$.
    \item $p\notin \pi$ and one of {\em(a)-(i)} below holds:
        \begin{enumerate}[{\em(a)}]

        \item $G=\PSL_n(q)$, $e(q,r)=r-1$, $(q^{r-1}-1)_r=r$,
            $\left[\frac{n}{r-1}\right]=\left[\frac{n}{r}\right]$ and for every $s\in \tau$ we have $e(q,s)=1$ and $n<s$. \label{L_epi-dpi_2B_a}
        \item $G=\PSU_n(q)$, $r\equiv1\,(\mathrm{mod}\, 4)$, $e(q,r)=r-1$, $(q^{r-1}-1)_r=r$, $\left[\frac{n}{r-1}\right]={\left[\frac{n}{r}\right]}$ and for every $s\in \tau$ we have $e(q,s)=2$ and $n<s$. \label{L_epi-dpi_2B_b}
        \item $G=\PSU_n(q)$, $r\equiv3\,(\mathrm{mod}\, 4)$, $e(q,r)=\frac{r-1}{2}$, $(q^{r-1}-1)_r=r$, $\left[\frac{n}{r-1}\right]={\left[\frac{n}{r}\right]}$ and for every $s\in \tau$ we have $e(q,s)=2$ and $n<s$. \label{L_epi-dpi_2B_c}
        \item $G=E_6(q)$, $\pi\cap\pi(G)\subseteq \pi(q-1)$, $3, 13 \in \pi\cap\pi(G)$, $5\notin \pi\cap\pi(G)$.\label{L_epi-dpi_2B_d}
        \item $G={}^2E_6(q)$, $\pi\cap\pi(G)\subseteq \pi(q+1)$, $3, 13 \in \pi\cap\pi(G)$, $5\notin \pi\cap\pi(G)$.
        \item $G=E_7(q)$, $\pi\cap\pi(G)$ is contained in one of the sets $\pi(q-1)$ or $\pi(q+1)$, $3, 13 \in \pi\cap\pi(G)$, $5,7 \notin \pi\cap\pi(G)$.
        \item $G=E_8(q)$, $\pi\cap\pi(G)$ is contained in one of the sets $\pi(q-1)$ or $\pi(q+1)$, $3, 13 \in \pi\cap\pi(G)$, $5,7 \notin \pi\cap\pi(G)$.
        \item $G=E_8(q)$, $\pi\cap\pi(G)$ is contained in one of the sets $\pi(q-1)$ or $\pi(q+1)$, $5, 31 \in \pi\cap\pi(G)$, $3,7 \notin \pi\cap\pi(G)$.
        \item $G=F_4(q)$, $\pi\cap\pi(G)$ is contained in one of the sets $\pi(q-1)$ or $\pi(q+1)$, $3, 13 \in \pi\cap\pi(G)$. \label{L_epi-dpi_2B_i}
        \end{enumerate}
    \end{enumerate}
\end{enumerate}
\end{lem}

In view of Lemmas \ref{max} and \ref{pinpi}, we need information about maximal subgroups of groups of Lie type. For information about maximal subgroups in classical groups we refer to \cite{KleiLie}.
Maximal subgroups of exceptional groups of Lie type are specified in Lemmas \ref{max_excep}, \ref{max_3D4} and \ref{max_2F4}.

Let $\overline{G}$~be an adjoint simple algebraic group of exceptional type $G_2$, $F_4$, $E_6$, $E_7$ or $E_8$ over $\overline{\F}_p$\,, the algebraic closure of the prime field $\F_p$\,, where $p$ is a prime. Let $\sigma$ be an endomorphism of $\overline{G}$ whose fixed point group $\overline{G}_\sigma=\{g\in \overline{G}\mid g^\sigma=g\}$ is finite. Then $\sigma$ is said to be a Frobenius morphism of~$\overline{G}$, and $G=(\overline{G}_\sigma)'$ is a finite simple exceptional group (exclude the cases $G_2(2)' \simeq \PSU_3(3)$ and ${}^2G_2(3)'\simeq \PSL_2(8)$). For a simple group of Lie type $N(q)$, let $rk(N(q))$ denote the untwisted Lie rank of $N(q)$ (i.e. the rank of the corresponding algebraic group).

\begin{lem}\label{max_excep} {\em \cite[Theorem 8]{LieSeitz}}
Let $N$~be a maximal subgroup of a finite exceptional group $\overline{G}_\sigma$ over $\F_q$, $q=p^n$. Then one of the following holds.
\begin{enumerate}[{\em(1)}]
\item $N=K_\sigma$ where $K$ is maximal closed $\sigma$-stable of positive dimension in $\overline{G}$; the possibilities are as follows:
    \begin{enumerate}[{\em(a)}]
    \item $K$~ $($and $N)$~is a parabolic subgroup.
    \item $K$~is reductive of maximal rank: the possibilities for $N$ are determined in {\em \cite{LieSaxlSeitz1992}}.
    \item $\overline{G}=E_7$, $p>2$ and $N=(2^2\times \P\Omega^+_8(q)\arbitraryext 2^2)\arbitraryext S_3$ or ${}^3D_4(q)\arbitraryext 3$.
    \item $\overline{G}=E_8$, $p>5$ and $N=\PGL_2(q)\times S_5$.
    \item $F^*(N)$ is as in  Table {\em \ref{Table3}}.
    \end{enumerate}
{\small
\begin{table}[htbp]
\begin{center}
\caption{ } \label{Table3}
\begin{tabular}{|c|l|}    \firsthline
$\overline{G}_\sigma'$ & \mbox{ possibilities for }$F^*(N)$\\ \hline
$G_2(q)$ & $A_1(q)$ $(p \ge 7)$\\ \hline
$F_4(q)$ & $A_1(q)$ ($p \ge 13$), $G_2(q)$ ($p=7$), $A_1(q)\times G_2(q)$ ($p\ge 3, q\ge 5$)\\ \hline
$E_6^\eta(q)$ & $A_2^\eta(q)$ ($p\ge 5$), $G_2(q)$ ($p\neq7$), $C_4(q)$ $(p\ge 3)$, $F_4(q)$,\\
                  & $A_2^\eta(q)\times G_2(q)$ ($( q,\eta) \neq(2,-) $)\\ \hline
$E_7(q)$ & $A_1(q)$ ($p \ge 17, 19$), $A_2^\eta(q)$ ($p\ge 5$), $A_1(q)\times A_1(q)$ $(p \ge 5)$, \\
                  & $A_1(q)\times G_2(q)$ ($p\ge 3, q\ge 5$), $A_1(q)\times F_4(q)$ $(q \ge 4)$, $G_2(q)\times C_3(q)$\\ \hline
$E_8(q)$ & $A_1(q)$ ($p \ge 23, 29, 31$), $B_2(q)$ ($p\ge 5$), $A_1(q)\times A_2^\eta(q)$ ($p\ge 5$), \\
 & $G_2(q)\times F_4(q)$, $A_1(q)\times G_2(q)\times G_2(q)$ ($p\ge 3, q\ge 5$), \\
  & $A_1(q)\times G_2(q^2)$ ($p\ge 3, q\ge 5$)\\ \hline
\end{tabular}
\end{center}
\end{table}
}
\item $N$ is of the same type as $\overline{G}$.
\item $N$~is an exotic local subgroup:
\[
    \begin{array}{llll}
    2^3\arbitraryext \SL_3(2)& < & G_2(p) & (p>2)\\
    3^3\arbitraryext \SL_3(3)& < & F_4(p) & (p \ge 5)\\
    3^{3+3}\arbitraryext \SL_3(3)& < & E_6^\eta(p)& (p \equiv \eta \pmod{3}, p\ge 5)\\
    5^{3}\arbitraryext \SL_3(5)& < & E_8(p^a) & (p\neq2, 5; a=1 \text{ or } 2, \text{ as } p^2 \equiv \pm 1 \pmod 5)\\
    2^{5+10}\arbitraryext \SL_5(2)& < & E_8(p)& (p >2).
    \end{array} \]
\item $\overline{G}=E_8$, $p>5$ and $N=(A_5\times A_6)\arbitraryext 2^2$.
\item $F^*(N)$~is simple and not a group of Lie type in characteristic $p$ {\em (}the possibilities for $F^*(N)$ are given up to isomorphism by {\em\cite{LieSeitz1999})}.
\item $F^*(N)=N(q_0)$~is a simple group of Lie type in characteristic $p$ with $rk(N(q_0))\le \frac{1}{2}rk(\overline{G})$ and one of the following holds:
    \begin{enumerate}[{\em(a)}]
    \item $q_0 \le 9$;
    \item $N(q_0)=A_2^\eta(16)$;
    \item $q_0 \le t(\overline{G})$ and $N(q_0)=A_1(q_0), {}^2B_2(q_0)$ or ${}^2G_2(q_0)$ {\em (}a constant $t(\overline{G})$ is defined in {\em \cite{LieSeitz1998})}.
    \end{enumerate} \label{L_max_excep_L(p)}
\end{enumerate}
\end{lem}

\begin{lem}\label{max_3D4} {\em \cite[Main theorem]{Kleidman1988}}
Let $G$ be a group with socle $L={}^3D_4(q)$, where $q=p^n$ and $p$ is prime. Assume that $M$ is a maximal subgroup of $G$ not containing $L$. Then $M_0=M\cap L$ is isomorphic to one of the following groups.
\begin{enumerate}[{\em(a)}]
\item $\ld q^{9}\rd\splitext (\SL_2(q^3)\circ (q-1))\arbitraryext d$, where $d=\gcd(2,q-1)$.
\item $\ld q^{11}\rd\splitext ((q^3-1)\circ \SL_2(q))\arbitraryext d$, where $d=\gcd(2,q-1)$.
\item $G_2(q)$.
\item $\PGL^\eta_3(q)$, $2<q\equiv \eta  \pmod{3}$.
\item ${}^3D_4(q_0)$, $q=q_0^\alpha$, $\alpha$ prime, $\alpha \neq 3$.
\item $\PSL_2(q^3)\times\PSL_2(q)$, $p=2$.
\item $(\SL_2(q^3)\circ\SL_2(q))\arbitraryext 2$, $p$ odd.
\item $((q^2+q+1)\circ \SL_3(q))\arbitraryext f_+ \arbitraryext 2$, where $f_+=\gcd(3,q^2+q+1)$.
\item $((q^2-q+1)\circ \SU_3(q))\arbitraryext f_- \arbitraryext 2$, where $f_-=\gcd(3,q^2+q+1)$.
\item $(q^2+q+1)^2\arbitraryext \SL_2(3)$.
\item $(q^2-q+1)^2\arbitraryext \SL_2(3)$.
\item $(q^4-q^2+1)\arbitraryext 4$.
\end{enumerate}
\end{lem}

We do not need a list of maximal subgroups of ${}^2B_2(q)$ because if $2\notin \pi$ then a $\pi$-Hall subgroup $H$ of ${}^2B_2(q)$ is abelian (see \cite[Lemma~14]{VdovinRevin2002}), and all maximal subgroups of ${}^2B_2(q)$ containing $H$ satisfy $\D_\pi$ by Lemma~\ref{L_Epi}(\ref{dpi_nil}). To simplify our proof of Theorem 1 we need a list of maximal subgroups of~${}^2F_4(q)$.

\begin{lem}\label{max_2F4}{\em \cite[Main Theorem]{Malle1991}}
Every maximal subgroup of $G={}^2F_4(q)$, $q=2^{2n+1}$, $n\geq1$, is isomorphic to one of the following.
\begin{enumerate}[{\em(a)}]
\item $\ld q^{11}\rd\splitext (A_1(q)\times (q-1))$.
\item $\ld q^{10}\rd\splitext ({}^2B_2(q)\times (q-1))$.
\item $\SU_3(q)\splitext 2$.
\item $((q+1)\times (q+1))\splitext \GL_2(3)$.
\item $((q-\sqrt{2q}+1)\times (q-\sqrt{2q}+1))\splitext \ld 96\rd$ if $q>8$.
\item $((q+\sqrt{2q}+1)\times (q+\sqrt{2q}+1))\splitext \ld 96\rd$.
\item $(q^2-\sqrt{2q}q+q-\sqrt{2q}+1)\splitext 12$.
\item $(q^2+\sqrt{2q}q+q+\sqrt{2q}+1)\splitext 12$.
\item $\PGU_3(q)\splitext 2$.
\item ${}^2B_2(q) \wr 2$.
\item $B_2(q)\splitext 2$.
\item ${}^2F_4(q_0)$, if $q_0=2^{2m+1}$ with $(2n+1)/(2m+1)$~prime.
\end{enumerate}
\end{lem}


We also need some information about automorphisms of groups of Lie type. Let $G$~be a simple group of Lie type. Definitions of diagonal, field and graph automorphisms of $G$ agrees with that of \cite{Steinberg}. The group of inner-diagonal automorphisms of $G$ is denoted by  $\widehat{G}$, while $\Aut(G)$ is the automorphism group of $G$.  By \cite[3.2]{Steinberg}, there exists a field automorphism $\rho$ of~$G$ such that every automorphism $\sigma$ of $G$ can be written $\sigma=\beta\rho^l\gamma$, with $\beta$ and  $\gamma$ being an inner-diagonal and a graph automorphisms, respectively, and $l\ge 0$. The group $\langle \rho\rangle$ is denoted by $\Phi_G$. In view of \cite[7-2]{GorLyons1983}, the group $\Phi_G$ is determined up to  $\widehat{G}$-conjugacy.  Since $G$ is centerless, we can identify $G$ with the group of its inner automorphisms.

\begin{lem}{\em  \cite[3.3, 3.4, 3.6]{Steinberg}}\label{automorphisms}
Assume that $G$~is a simple group of Lie type over $\F_q$ of characteristic $p$. Set $A=\Aut(G)$ and $\widehat{A}=\widehat{G}\Phi_G$. Then the following hold.
\begin{enumerate}[{\em(a)}]
\item $G\le \widehat{G} \le \widehat{A}\le A$ is a normal series for $A$.
\item  $\widehat{G}/G$ is abelian; $\widehat{G}=G$ for the groups $E_8(q), F_4(q), G_2(q),{}^3D_4(q)$, in other cases the order of $\widehat{G}/G$ is specified in Table~$\ref{d}$.
    \begin{table}[htbp]\center
\caption{} \label{d}
\begin{tabular}{|c|c|}  \hline
        $G$ & $|\widehat{G}/G|$ \\ \hline
        $A_l(q)$ & $\gcd(l+1,q-1)$ \\
        ${}^2A_l(q)$ & $\gcd(l+1,q+1)$ \\
        $B_l(q)$, $C_l(q)$, $E_7(q)$ & $\gcd(2,q-1)$ \\
        $D_l(q)$ & $\gcd(4,q^l-1)$ \\
        ${}^2D_l(q)$ & $\gcd(4,q^l+1)$ \\
        $E_6(q)$ & $\gcd(3,q-1)$ \\
        ${}^2E_6(q)$ & $\gcd(3,q+1)$ \\ \hline
        \end{tabular}
\end{table}
\item $A=\widehat{A}$ with the exceptions: $A/\widehat{A}$ has order $2$ if $G$ is $A_l(q)$ ($l\ge 2$),  $D_l(q)$  ($l\ge 5$) or $E_6(q)$, or if  $G$ is $B_2(q)$ or $F_4(q)$ and $q=2^{2n+1}$, or if  $G$ is $G_2(q)$ and $q=3^{2n+1}$; $A/\widehat{A}$ is isomorphic to $S_3$ if $G$ is $D_4(q)$.
\end{enumerate}
\end{lem}

\section{Proof of the main theorem}

In view of Lemma \ref{2inpi}, we may assume that $2\notin \pi$.  By Lemma \ref{max}, to prove the identity $\D_\pi=\U_\pi$, it is sufficient to prove that in each simple nonabelian $\D_\pi$-group $G$ all maximal subgroups, containing a $\pi$-Hall subgroup of $G$, satisfy~$\D_\pi$. This is true if $G$~is an alternating group, a sporadic simple group or a  simple group of Lie type in characteristic $p\in\pi$ because these groups satisfy~$\U_\pi$ by Lemma \ref{pinpi}. Thus we may assume that $G$~is a simple $\D_\pi$-group of Lie type with the base field  $\F_q$ of characteristic $p$ and  $p\not \in \pi$.

Throughout this section,  $H$~is a $\pi$-Hall subgroup of $G$, $M$~is a maximal subgroup of $G$ with  $H\le M$. Our goal is to show that $M$ satisfies $\D_\pi$. So by Lemma~\ref{L_Epi}(\ref{Dpi_comp_crit}) it is sufficient to prove that every nonabelian composition factor of~$M$ satisfies~$\D_\pi$. Since $M$ contains $H$, $M$ satisfies $\E_\pi$. Lemma~\ref{L_Epi}(\ref{Epi}) implies that every composition factor of $M$ satisfies  $\E_\pi$.

Let $r$ be the smallest prime in  ${\pi\cap \pi(G)}$ and $\tau=(\pi\cap\pi(G))\setminus\{r\}$.

\bigskip
{\bf Step 1. }If $S$ is a composition factor of $M$ and $S\in \E_\pi \setminus \D_\pi$, then $r\in \pi(S)$, $|\pi\cap\pi(S)|\ge 2 $ and $S\simeq \PSL^\eta_{n_1}(q_1)$.

\begin{proof}
Clearly, if $|\pi\cap\pi(S)|\le 1$, then $S\in \D_\pi$. Therefore we have that $|\pi\cap\pi(S)|\ge 2 $, in particular, $|\pi\cap\pi(G)|\ge 2 $.
 Since $G$ satisfies~$\D_\pi$, it follows from Lemma \ref{t.2.1} that $G$ satisfies $(*)$, i.e., every $\pi$-subgroup of $G$ has a normal abelian $\tau$-Hall subgroup. Lemma~\ref{L_Epi}(\ref{existpi}) implies that every  $\pi$-subgroup of~$S$ is a homomorphic image of a $\pi$-subgroup of~$G$. Hence a $\pi$-Hall subgroup of~$S$ has a normal abelian $\tau$-Hall subgroup.  If $r\notin \pi(S)$ then a $\pi$-Hall subgroup of~$S$ is abelian. So by Lemma~\ref{L_Epi}(\ref{dpi_nil}) we have $S\in \D_\pi$ and this is a contradiction with $S\in\E_\pi\setminus \D_\pi$. Therefore we conclude $r\in \pi(S)$, as required.  Thus every $\pi$-subgroup of~$S$ possesses a normal abelian $\tau$-Hall subgroup, i.e. $S$ satisfies~$(*)$.

 Since $S\in \E_\pi\setminus \D_\pi$, the possibilities for $S$ are determined in Lemma~\ref{epi-dpi_lie}. Suppose that $S$ satisfies item $(1)$ of Lemma~\ref{epi-dpi_lie}. Then $S\simeq O'N$ and $\pi\cap\pi(S)=\{3,5\}$.  But a $\{3,5\}$-Hall subgroup of $O'N$ does not possess a normal Sylow 5-subgroup, hence $O'N$ does not satisfy $(*)$ and this case is impossible.

  Suppose that $S$~is a group of Lie type with the base field $\F_{q_1}$ of characteristic~$p_1$. Assume first that $S$ satisfies item $2(A)$ of Lemma~\ref{epi-dpi_lie}, and so $p_1 \in \pi$. If $p_1\neq r$ then $r\in (\pi\cap\pi(S))\setminus\{p_1\}$ and $r$ does not divide $|W(S)|$. Since $\pi(|W(S)|)=\pi(l!)$ for some natural~$l$, we obtain that $l<r<p_1$ and it contradicts the fact that $p_1$ divides $|W(S)|$.  Suppose now that $p_1=r$. Denote by $U$ a Sylow $p_1$-subgroup of $S$. In view of \cite[Theorem~3.2]{Gross1986}, a Borel subgroup $B=N_S(U)$ contains a $\pi$-Hall subgroup $H_1$ of $S$. If $Q$ is a normal abelian $\tau$-Hall subgroup of $H_1$ then  $H_1=U\times Q$. Hence $H_1$ is nilpotent and by Lemma~\ref{L_Epi}(\ref{dpi_nil}) we have $S\in \D_\pi$, a contradiction with $S\in \E_\pi\setminus  \D_\pi$. Assume now that  $S$ satisfies item 2(B) of Lemma~\ref{epi-dpi_lie}, and so $p_1\notin \pi$. If $S$ satisfies one of items  \ref{L_epi-dpi_2B_d}-\ref{L_epi-dpi_2B_i}, then we have that $\pi\cap\pi(S)\subseteq \pi(q_1\pm 1)$ and therefore $e(q_1,t)=e(q_1,s)$ for every  $t,s \in \pi\cap\pi(S)$. Since $S$ satisfies $(*)$, $S$ satisfies $\D_\pi$ by Lemma~\ref{dpi_not2},  and this contradicts our assumption that $S\in \E_\pi\setminus  \D_\pi$. Thus $S$ satisfies one of items  \ref{L_epi-dpi_2B_a}-\ref{L_epi-dpi_2B_c} of Lemma~\ref{epi-dpi_lie}, and so $S\simeq \PSL_{n_1}^\eta(q_1)$.
\end{proof}

{\bf Step 2.}  If  $M$ is an almost simple group, then $M\in \D_\pi$.

\begin{proof}
  Let $S$ be nonabelian simple group with  $S\le M \le \Aut(S)$. As mentioned previously, $M\in \D_\pi$ if and only if $S\in \D_\pi$.  We derive a contradiction by assuming that $M\notin \D_\pi$ and thus $S\notin \D_\pi$.  By Step 1, we have $S\simeq \PSL_{n_1}^\eta(q_1)$ and $S$ satisfies one of items  \ref{L_epi-dpi_2B_a}-\ref{L_epi-dpi_2B_c} of Lemma~\ref{epi-dpi_lie}. Consider the structure of a $\pi$-Hall subgroup of $S$.
   Since $S\simeq \PSL_{n_1}^\eta(q_1)\in \E_\pi \setminus \D_\pi$, by Lemmas~\ref{L_Epi}(\ref{epi=cpi}) and~\ref{L_Epi}(\ref{Dpi_comp_crit}) we have that $\SL_{n_1}^{\eta}(q_1)$ and $\GL_{n_1}^{\eta}(q_1)$ lie in $\E_\pi\setminus \D_\pi$. By Lemma~\ref{epi-dpi_lie}, we have that $\gcd(n_1, q_1-\eta)_\pi=1$, therefore $|Z(\SL_{n_1}^{\eta}(q_1))|_{\pi}=1$, and so a $\pi$-Hall subgroup of $S$ is isomorphic to a $\pi$-Hall subgroup of $\SL_{n_1}^{\eta}(q_1)$.

    Now we specify the structure of a $\pi$-Hall subgroup of $\SL_{n_1}^\eta(q_1)$. We start with $\GL_{n_1}^\eta(q_1)$. Observe that $\pi(\GL_{n_1}^\eta(q_1))=\pi(S)$. Recall (by Step 1) that $r\in \pi(S)$. By Lemma~\ref{epi-dpi_lie}, for every $t\in \tau \cap \pi(S)$ if $\eta=+$ then $e(q_1,t)=1$,   and if $\eta=-$ then $e(q_1,t)=2$. Hence a $\tau$-Hall subgroup of $\GL_{n_1}^{\eta}(q_1)$ lies in a subgroup of diagonal matrices $D$, and $D$ is a direct product of $n_1$ copies of a cyclic group of order $q_1-\eta$. Since $\tau \cap \pi(S)\neq\varnothing$, it follows that $|q_1-\eta|\ge 5$, and so $N_{\GL_{n_1}^\eta(q_1)}(D)=D\splitext S_{n_1}$, where $S_{n_1}$ is a group of monomial matrices and $S_{n_1}$ acts on $D$ by permuting diagonal elements. Consider a Sylow $r$-subgroup~$R$ of~$S_{n_1}$. Lemma~\ref{epi-dpi_lie} implies that $r$ does not divide $q_1-\eta$, hence $R\le \SL_{n_1}^\eta(q_1)$.  Denote the number  $\left[\frac{n_1}{r}\right]$ by $k$. Since $\left[ \frac{n_1}{r}\right]=\left[ \frac{n_1}{r-1}\right]$, for some natural $d$ we have $n_1=kr+d=k(r-1)+(d+k)$ and $d+k<r-1$. Hence $k<r-1$ and $R\simeq r^k$. In view of Lemma~\ref{r4ast'}, we obtain that $|\SL_{n_1}^\eta(q_1))|_r=r^k$. Thus a $\pi$-Hall subgroup of $\SL_{n_1}^\eta(q_1)$, up to conjugation, lies in $(D\cap\SL_{n_1}^\eta(q_1))\splitext R$. If $D_1$ and $R_1$ are the images of $D\cap\SL_{n_1}^\eta(q_1)$ and $R$ under the canonical homomorphism $\SL_{n_1}^\eta(q_1)\rightarrow S$, respectively, then we may assume that a $\pi$-Hall subgroup of~$S$ lies in $$ D_1\splitext R_1\simeq (q_1-\eta)^{n_1-1}\splitext r^k.$$

    Recall that $M$ contains a $\pi$-Hall subgroup $H$ of $G$. So $H$ is a $\pi$-subgroup of  $\Aut(S)$, and since $\Aut(S)/\widehat{S}\Phi$ is a 2-group by Lemma~\ref{automorphisms}, where $\Phi$ denotes the group $\Phi_S$, we have that $H$ lies in $\widehat{S}\Phi$. Since $H$ is a $\pi$-Hall subgroup of~$M$, by Lemma~\ref{L_Epi}(\ref{Epi}) we have that $H\cap S$ is a $\pi$-Hall subgroup of $S$. Write $H_1=H\cap S$, and observe that by Lemma~\ref{L_Epi}(\ref{epi=cpi}) $S\in\C_\pi$, hence by Lemma~\ref{L_Epi}(\ref{prop4.8}) there exists  a $\pi$-Hall subgroup $H_0$ of $S\Phi$ with $H_0\cap S=H_1$. Now we want to show that it is possible to assume $H\le H_0$. By Lemma~\ref{automorphisms}, we have $|\widehat{S}\Phi:S\Phi|=\gcd(n_1,q_1-\eta)$ and it follows from $\gcd(n_1,q_1-\eta)_\pi=1$, that $H_0$ is a $\pi$-Hall subgroup of $\widehat{S}\Phi$.  Consider the normalizer $N_{\widehat{S}\Phi}(H_1)$. It is obvious that $N_{\widehat{S}\Phi}(H_1)$ contains both $H_0$ and $H$, and so $H_0$ is a $\pi$-Hall subgroup of $N_{\widehat{S}\Phi}(H_1)$. The normalizer $N_{\widehat{S}\Phi}(H_1)$ possesses a normal series $$H_1\unlhd N_{S}(H_1)\unlhd N_{\widehat{S}\Phi}(H_1),$$ where $N_S(H_1)/H_1$ is a $\pi'$-group and $N_{\widehat{S}\Phi}(H_1)/N_S(H_1)$ is solvable. Now Lemma~\ref{L_Epi}(\ref{Dpi_comp_crit}) implies that $N_{\widehat{S}\Phi}(H_1)$ satisfies $\D_\pi$. Since $H_0$ is a $\pi$-Hall subgroup of $N_{\widehat{S}\Phi}(H_1)$,  we may assume  $H\le H_0$, up to conjugation. In particular, $H=H_1\splitext \langle \varphi \rangle$, where $\langle \varphi \rangle$~is a $\pi$-subgroup of $\Phi$. Hence, by previous paragraph, $H$ is included in a subgroup $(D_1\splitext R_1)\splitext \langle \varphi \rangle$ of $S\splitext \langle \varphi \rangle$, and $(D_1\splitext R_1)\splitext \langle \varphi \rangle$ is isomorphic to $((q_1-\eta)^{n_1-1}\splitext r^k)\splitext \langle \varphi \rangle$. Moveover,   $\varphi$ centralizes $R_1$ and acts by  $x\mapsto x^{p_1^{\alpha}}$ on every group  $(q_1-\eta)$.   It follows from the action of $R_1$ on $D_1$, that  $|C_{D_1}(R_1)|_\pi\le (q_1-\eta)_{\pi}^{k+d-1}$. Thus $|C_{H}(R_1)|\le r^k(q_1-\eta)_{\pi}^{k+d-1}|\varphi|$.

On the other hand, in $\GL_{n_1}^\eta(q_1)$ there exists a subgroup $L$ of the following type:
\[\left(
\begin{array}{cc}
    \overbrace{ \begin{array}{ccc}
    \GL_{r-1}^\eta(q_1)& \cdots& 0\\
   \vdots &\ddots& \vdots\\
   0 &\cdots & \GL_{r-1}^\eta(q_1)
    \end{array}}^{k \text{ times}} & 0 \\
 0& \overbrace{ \begin{array}{ccc}
    q_1-\eta&\cdots &0\\
   \vdots &\ddots& \vdots\\
    0&\cdots & q_1-\eta
    \end{array}}^{k+d \text{ times}}
\end{array}
 \right)\]
Now $\GL_{r-1}^\eta(q_1)$ possesses a cyclic subgroup of order $q_1^{r-1}-1$. Since $r(q_1-\eta)$ divides $q_1^{r-1}-1$, we obtain that $\GL_{r-1}^\eta(q_1)$ has a cyclic subgroup of order $r(q_1-\eta)$. Thus $L$ (and hence $\GL_{n_1}^\eta(q_1)$) possesses a subgroup $r^k\times (q_1-\eta)^{2k+d}$. So in $L\cap\SL_{n_1}^\eta(q_1)$ there is an abelian subgroup of order $r^k(q_1-\eta)^{2k+d-1}$ and $L\cap\SL_{n_1}^\eta(q_1)$ is $\varphi$-invariant. It follows that $S\splitext \langle\varphi\rangle$ possesses a $\pi$-subgroup $F=(r^k\times (q_1-\eta)^{2k+d-1}_\pi)\splitext \langle\varphi\rangle$. Moreover, $|F|_r=|R_1|=|\SL_{n_1}^\eta(q_1)|_r$, hence up to conjugation we may assume that $R_1\le F$. Clearly $R_1$ lies in $Z(F)$ and $|F|>|C_{H}(R_1)|$. Thus there is a $\pi$-subgroup of $M$ (and of $G$) such that it is not isomorphic to any subgroup of a $\pi$-Hall subgroup $H$ of $G$, and this is the desired contradiction since $G\in \D_\pi$.
\end{proof}

Now we can assume that $M$ is not almost simple, and it seems quite natural to proceed our proof in two cases:  $G$ is a classical group and $G$ is an exceptional group.

\bigskip

{\bf Step 3.} If $G$ is a classical group and $M$ is not almost simple, then $M\in \D_\pi$.

\begin{proof}
  Since $G\in \D_\pi$, $|\pi\cap \pi(G)|\ge 2$ and $2, p\notin \pi$, we obtain that $G$ satisfies either Condition II or Condition~III of Lemma~\ref{Dpi}. If $G$ satisfies either item (g) or item (h) of Condition~II, then $H$ is cyclic and $M$ satisfies $\D_\pi$ by Lemma~\ref{L_Epi}(\ref{dpi_nil}). Suppose that $G$ satisfies either Condition~III or one of items (a)-(f) of Condition~II. Then $e(q,t)=e(q,s)$ for every $t,s\in \tau$. Set $$a=e(q,r) \mbox{ and } b=e(q,t) \mbox{ for every }t\in \tau.$$  Suppose that $M\notin \D_\pi$. Then by Lemma~\ref{L_Epi}(\ref{Dpi_comp_crit}) $M$ has a composition factor  $S$ which does not satisfy $\D_\pi$, and by Step 1, $S\simeq\PSL_{n_1}^\eta(q_1)$, $r\in \pi(S)$ and one of items  \ref{L_epi-dpi_2B_a}-\ref{L_epi-dpi_2B_c} of Lemma~\ref{epi-dpi_lie} holds for $S$. Set $$a_1=e(q_1,r) \mbox{ and } b_1=e(q_1,t) \mbox{ for every }t\in \tau \cap \pi(S).$$ Assume first that $q_1=q$. Since $S\in \E_\pi\setminus \D_\pi$,  we have that $a\neq b$. Hence  $G$ satisfies Condition II. Then $b\ge r>2$ and it contradicts the fact that $b\le 2$ which follows from items \ref{L_epi-dpi_2B_a}-\ref{L_epi-dpi_2B_c} of Lemma~\ref{epi-dpi_lie}.

  Assume now that $q_1\neq q$. Since $M$ is not almost simple, by famous Aschbacher's theorem~\cite{Asch}, we obtain that  $M$ belongs to one of Aschbacher's classes. According to Tables 3.5A--3.5F from~\cite{KleiLie}, it is sufficient to consider the following cases:

\begin{itemize}
\item[(a)] $G=\PSL_n(q)$, $M\simeq \left\ldbrack c \right\rdbrack  \arbitraryext \PSL_m(q^u) \arbitraryext \left\ldbrack d u \right\rdbrack $, where $n=mu$, $u$~is prime, $c=\frac{\gcd(q-1,m)(q^u-1)}{(q-1)\gcd(q-1,n)}$, ${d=\frac{\gcd(q^u-1,m)}{\gcd(q-1,m)}}$;
\item[(b)] $G=\PSL_n(q)$, $M\simeq \left\ldbrack \frac{c}{\gcd(q-1,n)} \right\rdbrack  \arbitraryext \PGL_n(q_0)$, where $q=q_0^u$, $u$~is prime, ${c=\frac{q-1}{\lcm\left( q_0-1,\frac{q-1}{\gcd(q-1,n)}\right)}}$;
\item[(c)] $G=\PSL_n(q)$, $M\simeq \PSU_n(q_0)\arbitraryext \left\ldbrack  \frac{\gcd(q_0+1,n)c}{\gcd(q-1,n)}\right\rdbrack$, where ${q=q_0^2}$, \\ ${c=\frac{q-1}{\lcm\left( q_0+1,\frac{q-1}{\gcd(q-1,n)}\right)}}$;
\item[(d)] $G=\PSU_n(q)$, $M\simeq  \left\ldbrack q^{m(2n-3m)} \splitext  c/\gcd(q+1,n) \right\rdbrack\arbitraryext(\PSL_m(q^2)\times \\\PSU_{n-2m}(q))\arbitraryext \left\ldbrack d \right\rdbrack$, where $1\le m \le [n/2]$, $$c=|\{(\lambda_1,\lambda_2)|\lambda_i\in \F_{q^2}, \lambda_2^{q+1}=1,\lambda_1^{m(q-1)}\lambda_2^{n-2m}=1\}|,$$ $$d=(q^2-1)\gcd(q^2-1,m)\gcd(q+1,n-2m)/c;$$
\item[(e)] $G=\PSU_n(q)$, $M\simeq \left\ldbrack  \frac{(q-1)\gcd(q+1,\frac{n}{2})}{\gcd(q+1,n)}\right\rdbrack \arbitraryext \PSL_{n/2}(q^2)\arbitraryext \left\ldbrack  \frac{\gcd(q^2-1,\frac{n}{2})}{\gcd(q+1,\frac{n}{2})} \right\rdbrack \arbitraryext 2 $;
\item[(f)] $G=\PSU_n(q)$, $M\simeq \left\ldbrack  c \right\rdbrack  \arbitraryext \PSU_m(q^u) \arbitraryext \left\ldbrack  d u \right\rdbrack $, where $n=mu$, $u$~is prime, $u\geq 3$, $c=\frac{\gcd(q+1,m)(q^u+1)}{(q+1)\gcd(q+1,n)}$, ${d=\frac{\gcd(q^u+1,m)}{\gcd(q+1,m)}}$;
\item[(g)] $G=\PSU_n(q)$, $M\simeq \left\ldbrack \frac{c}{\gcd(q+1,n)} \right\rdbrack \arbitraryext \PGU_n(q_0)$, where  ${q=q_0^u}$, $u$~is prime, $u\geq 3$, $c=\frac{q+1}{\lcm\left( q_0+1,\frac{q+1}{\gcd(q+1,n)}\right)}$.

\end{itemize}

In cases (a), (b), (d) and (e) $S$  is isomorphic to $\PSL_{n_1}(q_1)$ with $q_1\neq q$.  By Lemma~\ref{epi-dpi_lie} the following conditions hold:
\begin{itemize}
\item[$(1)$] $a_1=r-1$;
\item[$(2)$] $(q_1^{r-1}-1)_r=r$;
\item[$(3)$] $\left[ \frac{n_1}{r-1}\right]=\left[ \frac{n_1}{r}\right]$;
\item[$(4)$] $b_1=1$ and $n_1<t$ for every $t\in \tau \cap \pi(S)$.
\end{itemize}

It follows from $(1)$ and $(3)$ that $n_1\ge r$. Observe that if $k=\left[ \frac{n_1}{r}\right]=\left[ \frac{n_1}{r-1}\right]$ then $k<r-1$. Hence we have  $r\le n_1 < (r-1)^2$. By $(1)$ and $(4)$ we see that $(q_1-1)_r=1$ and $(q_1-1)_t\neq 1$ for every $t\in \tau \cap \pi(S)$.

\medskip
 {\it Case} (a): $S\simeq \PSL_m(q^u)$, where $mu=n$ and $u$~is prime. Note that $a_1=e(q^u,r)$, and so $a=e(q,r)\geq a_1>1$. By Lemma~\ref{r4ast'}, we have that  $b=e(q,t)$ divides $u$ and $(q_1-1)_t=(q^u-1)_t=(q^b-1)_t(\frac{u}{b})_t$. Hence $b$ is equal to 1 or~$u$. Assume  that  $G$ satisfies Condition III, and so $a=b$. Since $a>1$, we conclude that $a=u$. Then $a_1=1$, and it contradicts the fact that $a_1>1$. Thus we obtain that $G$ satisfies either item (a) or item (b) of Condition~II, and so $a=r-1$, $b=r=u$. Since $r=u$ divides $n$, we have that $G$ cannot satisfy item (b) of Condition II, where $n \equiv -1 \pmod{r}$. Hence $G$ satisfies item (a) of Condition II. Thus we have that $(q^{r-1}-1)_r=r$, $\left[ \frac{n}{r-1}\right]=\left[ \frac{n}{r}\right]=m$, and, in particular, $m<r-1$ and $\vert S\vert_r=1$. Using Lemma~\ref{r4ast'}, we can calculate $|G|_r$ and~$|M|_r$:  $$|G|_r=\frac{(q^{r-1}-1)_r^m(m!)_r}{(q-1)_r\gcd(q-1,n)_r}=r^m, |M|_r=r.$$
  As $m\ge 2$ we obtain $|G|_r\geq r^2$. Thus  $|G|_r > |M|_r$, which is a contradiction with the fact that $M$ contains a $\pi$-Hall subgroup of~$G$.

\medskip
 {\it Case} (b): $S\simeq \PSL_n(q_0)$, where $q=q_0^u$ and $u$~is prime. Note that $b_1=e(q_0,t)=1$, and so $b=e(q_0^u,t)=1$. Hence we have that  $G$ satisfies Condition III  and $a=b=1$.  Lemma~\ref{r4ast'} implies that $a_1=e(q_0,r)$ divides~$u$ and $(q-1)_r=(q_0^u-1)_r=(q_0^{a_1}-1)_r(\frac{u}{a_1})_r$. Since  $a_1$ is even ($a_1=r-1$) and $u$ is prime, we conclude that  $a_1=u=2$ and $r=3$. Using Lemma~\ref{r4ast'}, we can calculate $|G|_3$ and~$|M|_3$. $$|G|_3=\frac{(q-1)^{n-1}_3(n!)_3}{\gcd(q-1,n)_3}=\frac{(q_0^2-1)_3^{n-1}(n!)_3}{\gcd(q_0^2-1,n)_3};$$ $$|M|_3=\frac{(q-1)_3(q_0^2-1)^k_3(k!)_3}{\lcm\left( q_0-1,\frac{q-1}{\gcd(q-1,n)}\right)_3\gcd(q-1,n)_3(q_0-1)_3}, \text{ where } k=\left[\frac{n}{2}\right].$$ It follows from  $r\le n < (r-1)^2$ that $n=3$. Since $(q_0-1)_3=1$ and $(q_0^2-1)_3=3$,  we obtain that $|G|_3=3^2$  and $|M|_3=3$. Thus  $|G|_3 > |M|_3$, which is a contradiction with the fact that $M$ contains a $\pi$-Hall subgroup of~$G$.

\medskip
  {\it Cases} (d) and (e): $S\simeq \PSL_{n_1}(q^2)$. Note that $b_1=e(q^2,t)=1$, and so $b=e(q,t)$ is equal to 1 or 2.  Thus  $G$ satisfies Condition~III  and $a=b$. Then $a=e(q,r)$ equals 1 or 2, and so $a_1=e(q^2,r)=1$, which is a contradiction with the fact that $a_1=r-1\ge 2$.

In cases (c), (f) and (g) $S$ is isomorphic to $\PSU_{n_1}(q_1)$. By Lemma~\ref{epi-dpi_lie}  the following conditions hold:
\begin{itemize}
\item[$(1')$] either $r\equiv 1 \pmod{4}$ and $a_1=r-1$, or $r\equiv 3 \pmod{4}$ and $a_1=\frac{r-1}{2}$;
\item[$(2')$] $(q_1^{r-1}-1)_r=r$;
\item[$(3')$] $\left[ \frac{n_1}{r-1}\right]=\left[ \frac{n_1}{r}\right]$;
\item[$(4')$] $b_1=2$ and $n_1<t$ for every $t\in \tau \cap \pi(S)$.
\end{itemize}

In view of $(1')$, we see that either $a_1\equiv 0 \pmod{4}$, or $a_1\equiv 1 \pmod{2}$. So $a_1$ cannot equal $2k$ with $k$~odd, in particular $a_1\neq 2$. It is obvious that in this case the inequalities $r\le n_1 < (r-1)^2$ are also valid for $n_1$. By  $(4')$ we have that ${(q_1-1)_t=1}$ and $(q_1^2-1)_t\neq 1$ for every $t\in \tau \cap \pi(S)$.

\medskip
{\it Case} (c): $S\simeq\PSU_n(q_0)$ and  $q=q_0^2$. Note that $b_1=e(q_0,t)=2$, and so $b=e(q_o^2,t)=1$. Hence $G$ satisfies Condition III,  and so $a=b=1$. It follows from $a=e(q_0^2,r)=1$, that $a_1=e(q_0,r)$ equals 1 or 2.  As we mentioned before, $a_1$ cannot equal 2, and thus $a_1=1$ and $r=3$. Using Lemma~\ref{r4ast'}, we can calculate $|G|_3$ and~$|M|_3$. As in case (b) we have $n=3$ and $|G|=3^2$. Since $a_1=1$ we obtain that $a_1^*=2$. $$|M|_3=\frac{(q_0^2-1)_3(q-1)_3}{(q_0+1)_3\gcd(q-1,3)_3\lcm\left( q_0+1,\frac{q-1}{\gcd(q-1,3)}\right)_3}$$ It follows from  $(1')$ and $(2')$ that $(q_0-1)_3=(q_0^2-1)_3=3$, and so $(q_0+1)_3=1$. Hence we have that $|M|_3=3$. Thus  $|G|_3 > |M|_3$, and it is a contradiction with the fact that $M$ contains a $\pi$-Hall subgroup  of~$G$.

\medskip
 {\it Case} (f):  $S\simeq\PSU_m(q^u)$, where $mu=n$, $u$~is prime and $u\geq 3$. Note that $b_1=e(q^u,t)$, and so $b=e(q,t)\ge b_1=2$ and $b\neq u$.  By Lemma~\ref{r4ast'}, we have that $b$ divides $2u$ and $(q^{2u}-1)_t=(q^b-1)_t(\frac{2u}{b})_t$. Thus we have that $b$ equals 2 or $2u$. Assume that $G$ satisfies Condition~III, and so $a=b$.  Since $a_1=e(q^u,r)$ cannot equal 2, we have that $a=e(q,r)\neq 2u$, and thus $a=b=2$. As $a\geq a_1$, we have $r=3$ and $a_1=1$.  Since $(q^u-1)_3\neq 1$, it follows from Lemma~\ref{r4ast'} that $a$ divides $u$ and  $(q^{u}-1)_3=(q^2-1)_3(\frac{u}{2})_3$, and this is a contradiction with the fact that $u$ is an odd prime.  Thus we obtain that $G$ satisfies one of items (c)-(f) of Condition~II, and so  $b=2r=2u$. Since $r=u$ divides $n$, we have that $G$ cannot satisfy items (e) and (f) of Condition~II,  where $n\equiv -1 \pmod r$. Hence $G$ satisfies either item (c) or item (d) of Condition II. Thus we have that $a^*=a_1^*=r-1$, $(q^{r-1}-1)_r=r$,  $\left[\frac{n}{r-1}\right]=\left[\frac{n}{r}\right]=m$,  in particular $m<r-1$ and $\vert S\vert_r=1$. As in the case~(a) we have  $|G|_r\geq r^2$ and $|M|_r=r$. Thus  $|G|_r > |M|_r$, which is a contradiction with the fact that $M$ contains a $\pi$-Hall subgroup of~$G$.

\medskip
  {\it Case} (g): $S\simeq \PSU_n(q_0)$, where $q=q_0^u$, $u$ is prime and $u\ge 3$. Note that $b_1=e(q_0,t)$ and $b_1\ge b=e(q_0^u,t)$. By Lemma~\ref{r4ast'}, we have that   $b_1=2$ divides~$ub$ and  $(q^b-1)_t=(q_0^{ub}-1)_t=(q_0^{2}-1)_t(\frac{ub}{2})_t$. Since  $u$~is an odd prime, we conclude that $b=2$. Thus $G$ satisfies Condition III  and $a=b=2$. Since $a=e(q_0^u,r)=2$, we have that $a_1=e(q_0,r)\neq u$ and $a_1\ge a$. Lemma~\ref{r4ast'} implies that  $a_1$ divides $2u$ and $(q^2-1)_r=(q_0^{2u}-1)_r=(q_0^{a_1}-1)_r(\frac{2u}{a_1})_r$. Then $a_1$ equals 2 or  $2u$, and it is impossible: as we noted earlier $a_1$ cannot be equal to $2k$ with $k$ odd.

  Thus we prove that in all cases (a)-(g) a composition factor $S\simeq \PSL_{n_1}^\eta(q_1)$ cannot satisfy statements \ref{L_epi-dpi_2B_a}-\ref{L_epi-dpi_2B_c} of Lemma~\ref{epi-dpi_lie}, and so $S\in \D_\pi$.
\end{proof}

{\bf Step 4.} If $G$ is an exceptional group, not isomorphic to ${}^2B_2(q)$, ${}^2F_4(q)$, ${}^2G_2(q)$, and $M$ is not almost simple, then $M\in \D_\pi$.

\begin{proof}
Since $G$ is not isomorphic to ${}^2B_2(q)$, ${}^2F_4(q)$, ${}^2G_2(q)$, by \cite[Lemmas 7-13]{VdovinRevin2002} we have that $H$ is abelian or $\pi\cap\pi(G)\subseteq \pi(q\pm1)$. If $H$ is abelian then $M$ satisfies $\D_\pi$ by Lemma~\ref{L_Epi}(\ref{dpi_nil}). Assume  now that  $\pi\cap\pi(G)\subseteq \pi(q\pm1)$. Suppose that $M\notin \D_\pi$. In view of Lemma~\ref{L_Epi}(\ref{Dpi_comp_crit}) and Step 1, $M$ possesses a composition factor $S\simeq \PSL_{n_1}^\eta(q_1)\in \E_\pi\setminus \D_\pi$. If $q_1=q^u$ for some natural $u$, then $\pi\cap\pi(S)\subseteq\pi(q\pm 1)\subseteq\pi(q^u\pm 1)$. It follows that $e(q^u,r)=e(q^u,t)$ for every $t\in \tau \cap \pi(S)$. By proof of Step~1, we have that $S$ satisfies $(*)$. Thus Lemma~\ref{dpi_not2} implies that  $S\in\D_\pi$, and it is a contradiction with our assumption that $S\in \E_\pi\setminus \D_\pi$. Assume now that $q_1\neq q^u$. It follows from Lemma~\ref{max_3D4} that all maximal subgroups of $G={}^3D_4(q)$ do not possess a composition factor isomorphic  $\PSL_{n_1}^\eta(q_1)$ with $q_1\neq q^u$.  If $G\neq {}^3D_4(q)$, then in Lemma~\ref{max_excep} we can find maximal subgroups of a finite exceptional group $\overline{G}_\sigma$ with $G=(\overline{G}_\sigma)'$. It is clearly that $M=N\cap G$, where $N$ is a maximal subgroup of $\overline{G}_\sigma$ not containing~$G$. Observe that $M\unlhd N$ and $S$ is a composition factor of $N$.  By Lemma~\ref{max_excep}, all possibilities for $N$, which possesses a composition factor $S$ isomorphic $\PSL_{n_1}^\eta(q_1)$ with $q_1\neq q^u$, are the following.

\begin{enumerate}[(a)]
 \item the possibilities for $N$~are determined in \cite{LieSaxlSeitz1992};
  \item $N$~is an exotic local subgroup;
 \item $F^*(N)$~is  simple.
 \end{enumerate}

 If $N$ satisfies the first statement, then by \cite{LieSaxlSeitz1992} we have that $S$ is isomorphic to $\PSL_2(5)$, $\PSL_3(2)$ or $\PSU_4(2)$. Since there exists no odd prime $t>3$ such that $q_1 \equiv 1 \pmod{t}$ or $2^2 \equiv 1 \pmod{t}$, we conclude that $S$ does not satisfy items \ref{L_epi-dpi_2B_a}-\ref{L_epi-dpi_2B_c} of Lemma~\ref{epi-dpi_lie}. Then $S\in\D_\pi$, and this is a contradiction with our assumption that $S\in \E_\pi\setminus \D_\pi$.

    If $N$ satisfies the second statement, then $S$ is isomorphic to one of the groups $\PSL_3(2)$, $\PSL_3(3)$, $\PSL_3(5)$ or $\PSL_5(2)$. Thus in this case we also have that there is no odd prime $t$ with $q_1 \equiv 1 \pmod{t}$. Hence $S$ cannot satisfy \ref{L_epi-dpi_2B_a} of Lemma~\ref{epi-dpi_lie}, and so $S\in \D_\pi$, which is a contradiction with $S\in \E_\pi\setminus \D_\pi$.

     Finally, if $F^*(N)$ is simple, then $F^*(M)$ is simple, a contradiction with our  assumption that $M$ is not almost simple.

     Thus we prove that in all cases (a)-(c) a composition factor $S\simeq \PSL_{n_1}^\eta(q_1)$ cannot satisfy statements \ref{L_epi-dpi_2B_a}-\ref{L_epi-dpi_2B_c} of Lemma~\ref{epi-dpi_lie}, and so $S\in \D_\pi$.
\end{proof}

    {\bf Step 5.} If $G$ is one of the groups ${}^2B_2(q)$, ${}^2G_2(q)$ or ${}^2F_4(q)$, then $M\in \D_\pi$.

 \begin{proof}
 The structure of a $\pi$-Hall subgroup $H$ of $G$ was obtained in \cite{VdovinRevin2002}. By  \cite[Lemma 14]{VdovinRevin2002}, if $G$ is either ${}^2B_2(q)$ or ${}^2G_2(q)$, or if $G$ is ${}^2F_4(q)$ and  $3\notin \pi$, then $H$ is abelian. Then it follows by Lemma~\ref{L_Epi}(\ref{dpi_nil}) that $M\in \D_\pi$.

  Assume now that $3\in\pi$ and $G$ is ${}^2F_4(q)$, where $q=2^{2m+1}$.  So we obtain that $r=3$.  Since $G\in \D_\pi$ and $|\pi\cap\pi(G)|\ge 2$, we have that $G\neq {}^2F_4(2)'$ by Lemma~\ref{Dpi}. All maximal subgroups of ${}^2F_4(q)$ are specified in Lemma~\ref{max_2F4}. Suppose that $M\notin \D_\pi$. In view of Step 1, we should consider maximal subgroups of ${}^2F_4(q)$, which have a composition factor $S$ isomorphic to $\PSL_{n_1}^\eta(q_1)$. It follows from Lemma~\ref{max_2F4} that~$S$ is isomorphic to $\PSL_2(q)$ or $\PSU_3(q)$. If $S\simeq \PSL_2(q)$ then $S$ cannot satisfy  \ref{L_epi-dpi_2B_a} of Lemma~\ref{epi-dpi_lie} because for $n_1=2$ and $r=3$ the inequality $ n_1\ge r$ does not hold. Hence $S\simeq \PSL_2(q)$ satisfies~$\D_\pi$. If $S\simeq \PSU_3(q)$ then $S$ also satisfies $\D_\pi$ because $2^{2m+1}\not \equiv 1\pmod{3}$ and $S$ cannot satisfy \ref{L_epi-dpi_2B_c} of Lemma~\ref{epi-dpi_lie}. So we prove that all composition factors of $M$  satisfy~$\D_\pi$.
 \end{proof}

     Thus we prove if $G$ is a simple $\D_\pi$-group of Lie type over $\F_q$ of characteristic $p$ and $2,p\notin \pi$, then every maximal subgroup $M$ of $G$ containing a $\pi$-Hall subgroup $H$ of $G$ satisfies $\D_\pi$. So we have $\D_\pi=\U_\pi$ by Lemma~\ref{max}.

\section*{Acknowledgments}

The author is grateful to E.P. Vdovin and D.O. Revin for stimulating and helpful discussions.

\end{document}